\documentclass[11pt,a4paper]{article}
\usepackage{graphicx}
\usepackage{amsmath}
\usepackage{amssymb}
\usepackage[width=16cm,height=24cm]{geometry}
\usepackage{cite}
\usepackage{subfigure}

\newcommand{\vect}[1]{\mathbf{#1}}
\newcommand{\mat}[1]{\mathbf{#1}}

\newcommand{\E}[1]{{\mathbb{E}}\{#1\}}

\newcommand{\etr}[1]{\operatorname{etr}(#1)}
\newcommand{\tr}[1]{\operatorname{tr}(#1)}
\newcommand{\diag}[1]{\operatorname{diag}(#1)}
\newcommand{\cond}[1]{\operatorname{cond}(#1)}
\renewcommand{\det}[1]{|#1|}
\newcommand{\loss}[1]{\mathsf{L}(#1)}
\newcommand{\losslimit}[1]{\mathsf{L}_{\infty}(#1)}
\newcommand{\risk}[1]{\mathsf{R}(#1)}
\newcommand{\sqnorm}[1]{\left\|#1\right\|^{2}}
\newcommand{\half}{\frac12}
\newcommand{\lambdamax}{\lambda_{\text{\tiny{max}}}}

\newcommand{\vg}{\vect{g}}

\newcommand{\vx}{\vect{x}}

\newcommand{\mD}{\mat{D}}
\newcommand{\mG}{\mat{G}}
\newcommand{\mGtilde}{\tilde{\mG}}
\newcommand{\mGbar}{\bar{\mG}}
\newcommand{\mH}{\mat{H}}
\newcommand{\mHtilde}{\tilde{\mH}}
\newcommand{\eye}[1]{\mat{I}_{#1}}
\newcommand{\mL}{\mat{L}}
\newcommand{\mQ}{\mat{Q}}
\newcommand{\mS}{\mat{S}}
\newcommand{\mU}{\mat{U}}
\newcommand{\mW}{\mat{W}}
\newcommand{\mWbar}{\bar{\mW}}
\newcommand{\mX}{\mat{X}}
\newcommand{\mXbar}{\bar{\mX}}
\newcommand{\mY}{\mat{Y}}
\newcommand{\mzero}{\mat{0}}

\newcommand{\mLambda}{\mat{\Lambda}}
\newcommand{\mPi}{\mat{\Pi}}
\newcommand{\mSigma}{\mat{\Sigma}}

\newcommand{\mSigmahat}{\hat{\mSigma}}

\newcommand{\Gaussian}[2]{\mathcal{N}\left(#1,#2\right)}

\newcommand{\mCN}[4]{\mathcal{N}_{#1}\left(#2,#3,#4\right)}

\newcommand{\chisquare}[1]{\chi^{2}_{#1}}
\newcommand{\dist}{\overset{d}{=}}

\begin{document}
\title{Covariance matrix estimation in the singular case using regularized Cholesky factor}
\author{Olivier Besson\thanks{ISAE-SUPAERO, Universit\'{e} de Toulouse, 10 avenue Marc P\'{e}legrin, 31055 Toulouse, France. Email: olivier.besson@isae-supaero.fr}}
\maketitle
\begin{abstract}
We consider estimating the population covariance matrix when the number of available samples is less than the size of the observations. The sample covariance matrix (SCM) being singular, regularization is mandatory in this case. For this purpose we consider minimizing Stein's loss function and we investigate a method based on augmenting the partial Cholesky decomposition of the SCM. We first derive the finite sample optimum estimator which minimizes the loss for each data realization, then the Oracle estimator which minimizes the risk, i.e., the average value of the loss. Finally a practical scheme is presented where the missing part of the Cholesky decomposition is filled. We conduct a numerical performance study of the proposed method and compare it with available related methods. In particular we investigate the influence of the condition number of the covariance matrix as well as of the shape of its spectrum.
\end{abstract}
\section{Introduction}
The need for accurate covariance matrix estimation, whether for multivariate analysis, data decorrelation or adaptive filtering is pregnant in many fields such as economics or engineering \cite{Pourahmadi13,Zagidullina21,Ledoit22b}. When the number of samples $n$ available to estimate the population covariance matrix $\mSigma$ is slightly above the size of the observations $p$, the sample covariance matrix (SCM) is known to become unreliable, and accordingly it is advisable to  regularize it. In the singular case, i.e., when $ n<p$, the case we consider herein, the SCM is singular and hence no longer positive definite as should a covariance matrix be, in which case regularization becomes mandatory. 

The most widespread method consists in linear shrinkage of the SCM, i.e,. a weighted linear combination of the SCM and  a target matrix, typically the identity matrix (the latter technique is often referred to as diagonal loading in the engineering community). When the target matrix is the identity matrix, linear shrinkage actually belongs to the class of orthogonally invariant estimators (OIE) which retain the SCM eigenvectors and apply a transformation to its eigenvalues: in this case linear shrinkage corresponds to an affine transformation of the SCM eigenvalues.  More precisely, let $\mS$ be the SCM and let its eigenvalue decomposition be $\mS=\mU\mLambda\mU^{T}$ with $\mU$ the matrix of eigenvectors and $\mLambda=\diag{\lambda_{1},\cdots,\lambda_{p}}$ the diagonal matrix of eigenvalues. OIE are of the form $\mSigmahat=\mU\mD\mU^{T}$ with $\mD=\diag{d_{1},\cdots,d_{p}}$ a diagonal matrix constructed from $\mLambda$. OIE have been extensively studied, see e.g., \cite{Stein86,Dey85,Dey86,Sheena92,Perron92} for early works. More recently Ledoit and Wolf  in a series of papers \cite{Ledoit04,Ledoit12,Ledoit17,Ledoit18,Ledoit20,Ledoit22} have thoroughly investigated them under the statistical framework  of random matrix theory (RMT), i.e., $p,n \rightarrow \infty$ with $p/n \rightarrow c$. Their methods mostly rely on RMT coupled with Stein's approach \cite{Stein56,James61,Stein81,Stein86}. Briefly stated, an estimate of a specific form (in their case $\mSigmahat=\mU\mD\mU^{T}$) is chosen along with a loss function $\loss{\mSigmahat,\mSigma}$. Ledoit and Wolf define what they refer to as finite sample optimal (FSOPT) the estimate obtained by minimizing $\loss{\mSigmahat,\mSigma}$: obviously this estimate is hypothetical since it depends on $\mSigma$. Tables 1, 2 and 3 of \cite{Ledoit21} provide the expressions of the FSOPT estimator for a large number of loss functions. The next step is to derive the asymptotic -as $p,n \rightarrow \infty$ with $p/n \rightarrow c$- limit of $\loss{\mSigmahat,\mSigma}$, say $\losslimit{\mSigmahat,\mSigma}$. Then the Oracle estimator is defined as the one which minimizes this limit. Note that in a classical framework where asymptotic is understood as $p$ fixed and $n \rightarrow \infty$, $\losslimit{\mSigmahat,\mSigma}$ is usually replaced by $\risk{\mSigmahat,\mSigma}=\E{\loss{\mSigmahat,\mSigma}}$. Finally, to come up with a practical scheme, a consistent (unbiased in the classical framework) estimate of $\losslimit{\mSigmahat,\mSigma}$ is derived, whose minimization leads to the final estimator. Proceeding along these lines, reference \cite{Ledoit04} derives a linear-shrinkage-based  bona fide estimator which now constitutes a reference estimator, see also \cite{Bodnar14} for a similar approach and a close estimator. While linear shrinkage (LS) is simple and quite effective, it can be improved upon using more complex transformations such as nonlinear shrinkage (NLS)   \cite{Ledoit12,Ledoit17,Ledoit18,Ledoit20} or possibly quadratic inverse shrinkage \cite{Ledoit22}.

All above referred techniques rely on keeping the eigenvectors of the SCM and modifying its eigenvalues. Another solution consists in regularizing the Cholesky factor of the SCM \cite{James61}, that is if $\mS = \mG \mG^{T}$ denotes the Cholesky decomposition of $\mS$ where $\mG$ is lower triangular with positive diagonal elements, the estimates are of the form $\mSigmahat=\mG \mD \mG^{T}$ with $\mD$ a diagonal matrix. While deemed less performant than OIE this type of estimators has two merits. First they are usually simpler, as Cholesky factors can be computed more efficiently than eigenvectors and eigenvalues. Besides, this approach enables one to obtain analytical expressions of the regularization parameters for a certain number of risks, see \cite{James61,Selliah64,Eaton87} due to the fact that the risk $\E{\loss{\mSigmahat,\mSigma}}$ no longer depends on unknown parameters, yielding closed-form expressions for $\mD$. As a result, the risk itself is minimized, not just an unbiased estimate of the risk. The above mentioned approaches deal with the case $n \geq p$  but the paper \cite{Tsukuma16c} considers the  case $n <p$. However, the regularized SCM has rank $n$ and is thus singular, which is problematic for an estimate of a covariance matrix which, in principle, is positive definite. Another class of regularization schemes is based on the modified Cholesky decomposition and relies on the fact that the elements of the Cholesky factor can be obtained from the coefficients of a sequence of  regression models \cite{Pourahmadi99,Pourahmadi00,Wu03}. Based on this property, reference \cite{Huang06} considers a penalized likelihood approach to regularize the Cholesky decomposition of the inverse of the covariance matrix. However this method is not applicable when $n < p$. In order to cope with the singular case, banding the Cholesky factor is proposed in \cite{Bickel08,Rothman10}. More precisely a sequence of reduced-dimension linear regression models is computed for the Cholesky factor of the covariance matrix \cite{Rothman10} or of its inverse \cite{Bickel08}. These methods work in the singular case and provide a positive definite covariance matrix estimate. Herein, we also aim at deriving a full rank Cholesky-based covariance matrix estimate in the singular case. Somehow information is missing when $n<p$, actually the $(p-n) \times  (p-n)$ right lower tail of the Cholesky factor of the population covariance matrix is not identifiable. Therefore, there is a need to fill this unknown part and a method is proposed to fulfill this need.

\section{Covariance matrix estimation}
\subsection{Outline}
Let $\mX=\begin{bmatrix} \vx_{1} & \vx_{2} & \cdots & \vx_{n} \end{bmatrix}$ be the $p \times n$ data matrix whose columns $\vx_{k}$ are independent and drawn from a multivariate Gaussian distribution with zero mean and full-rank population covariance matrix $\mSigma$. We assume here that $n <p$ so that $\mS=\mX\mX^{T}$ has rank $n$ with probability one. We let
\begin{equation}
\mS = \mG \mG^{T} = \begin{pmatrix} \mG_{11} \\ \mG_{21} \end{pmatrix}  \begin{pmatrix} \mG_{11}^{T} & \mG_{21}^{T} \end{pmatrix} 
\end{equation}
denote the Cholesky decomposition of $\mS$ where $\mG_{11}$ is a $n \times n$ lower triangular matrix with positive diagonal entries. Below we consider estimates of the form
\begin{equation}
\mSigmahat = \mGtilde \mD \mGtilde^{T} = \begin{pmatrix} \mG_{11} & \mzero  \\ \mG_{21} & \mGtilde_{22} \end{pmatrix} \begin{pmatrix} \mD_{1} & \mzero \\ \mzero & \mD_{2} \end{pmatrix}  \begin{pmatrix} \mG_{11}^{T} &  \mG_{21}^{T} \\ \mzero & \mGtilde_{22}^{T} \end{pmatrix}
\end{equation}
where $\mD=\diag{d_{1},\ldots,d_{n},d_{n+1},\ldots,d_{p}}$ is diagonal and $\mGtilde_{22}$ is a lower triangular matrix with positive diagonal entries. We wish to examine such estimates under Stein's loss
\begin{equation}\label{loss_function}
\loss{\mSigmahat,\mSigma}  = \tr{\mSigmahat\mSigma^{-1}} - \log \det{\mSigmahat\mSigma^{-1}}
\end{equation}
Towards this end we define the Cholesky decomposition of $\mSigma$  as
\begin{equation}
\mSigma = \mL \mL^{T} = \begin{pmatrix} \mL_{11} & \mzero \\  \mL_{21}  & \mL_{22} \end{pmatrix}  \begin{pmatrix} \mL_{11}^{T} & \mL_{21}^{T}  \\ \mzero & \mL_{22}^{T}\end{pmatrix}
\end{equation}
where $\mL_{11}$ and $\mL_{22}$ are $n \times n$ and $(p-n) \times (p-n)$ lower triangular matrices with positive diagonal entries. Note that due to $n <p$ $\mL_{22}$ is not identifiable from $\mS$. Actually, if we decompose $\mSigma$ as
\begin{equation}
\mSigma = \begin{pmatrix} \mSigma_{11} & \mSigma_{12} \\ \mSigma_{21} & \mSigma_{22} \end{pmatrix}  = \begin{pmatrix} \mL_{11}\mL_{11} ^{T} & \mL_{11} \mL_{21}^{T} \\ \mL_{21}\mL_{11} ^{T} & \mL_{21}\mL_{21}^{T}+\mL_{22}\mL_{22}^{T}  \end{pmatrix}
\end{equation}
and define $\mSigma_{2.1}=\mSigma_{22}-\mSigma_{21}\mSigma_{11}^{-1}\mSigma_{12}=\mL_{22}\mL_{22}^{T}$, the likelihood of the observations is
\begin{align}
L(\mX;\mSigma) &\propto \det{\mSigma}^{-\frac n2} \etr{-\half\mSigma^{-1}\mS} \nonumber \\
&= \det{\mSigma_{11}}^{-\frac n2} \det{\mSigma_{2.1}}^{-\frac n2}\etr{-\half \sqnorm{\mL^{-1}\mG}} \nonumber \\
&= \det{\mL_{11}}^{-n} \det{\mL_{22}}^{-n}\etr{-\half \sqnorm{\mL_{11}^{-1}\mG_{11}}}  \nonumber \\
&\times \etr{-\half \sqnorm{\mL_{22}^{-1}(\mG_{21}-\mL_{21}\mL_{11}^{-1}\mG_{11})}}
\end{align}
where $\etr{.}$ stands for the exponential of the trace of a matrix. Clearly the maximum is achieved at $\mL_{11}=n^{-1/2}\mG_{11}$, $\mL_{21}=n^{-1/2}\mG_{21}$ and
\begin{equation}
\underset{\mL_{11},\mL_{21}}{\max} L(\mX;\mSigma) \propto \det{\mL_{22}}^{-n}
\end{equation}
which is unbounded, revealing that $\mL_{22}$ is not identifiable. Consequently, nothing can be inferred about $\mL_{22}$ from the data and therefore we will need to fill in the missing information. In the sequel we will successively derive the FSOPT estimator which minimizes $\loss{\mSigmahat,\mSigma}$, then the Oracle estimator which minimizes $\E{\loss{\mSigmahat,\mSigma}}$ before proposing  a bona fide estimator that can be implemented in practice.
\subsection{Finite-sample optimum estimator}
The FSOPT estimator aims at minimizing $\loss{\mGtilde \mD \mGtilde^{T},\mSigma}$ with respect to $\mD$ and $\mGtilde_{22}$. First we rewrite the loss function in a more convenient form. We have
\begin{align}
\mL^{-1}\mGtilde &= \begin{pmatrix} \mL_{11}^{-1} & \mzero \\ -\mL_{22}^{-1}\mL_{21}\mL_{11}^{-1} & \mL_{22}^{-1} \end{pmatrix} \begin{pmatrix} \mG_{11} & \mzero  \\ \mG_{21} & \mGtilde_{22} \end{pmatrix} \nonumber \\
&= \begin{pmatrix} \mL_{11}^{-1} \mG_{11}& \mzero \\ \mL_{22}^{-1}(\mG_{21}-\mL_{21}\mL_{11}^{-1}\mG_{11}) & \mL_{22}^{-1} \mGtilde_{22}\end{pmatrix}
\end{align}
so that the loss function can be written as
\begin{align}\label{loss_function}
\loss{\mGtilde \mD \mGtilde^{T},\mSigma} &= \tr{\mL^{-1}\mGtilde \mD \mGtilde^{T}\mL^{-T}} - \log \det{\mL^{-1}\mGtilde \mD \mGtilde^{T}\mL^{-T}} \nonumber \\
&= \tr{\mL^{-1}\mG \mD_{1} \mG^{T}\mL^{-T}} - \log \det{\mD_{1} } - \log \det{\mL_{11}^{-1}\mG_{11} \mG_{11}^{T}\mL_{11}^{-T}} \nonumber \\
&+ \tr{\mL_{22}^{-1}\mGtilde_{22} \mD_{2} \mGtilde_{22}^{T}\mL_{22}^{-T}} - \log \det{\mL_{22}^{-1}\mGtilde_{22} \mD_{2} \mGtilde_{22}^{T}\mL_{22}^{-T}}
\end{align}
Observing that 
\begin{equation}
\tr{\mL^{-1}\mG \mD_{1} \mG^{T}\mL^{-T}} = \sum_{j=1}^{n} d_{j} \vg_{j}^{T}\mSigma^{-1}\vg_{j}
\end{equation}
where $\vg_{j}$ is the $j^{\text{th}}$ column of $\mG$, we find that minimization with respect to $\mD_{1}$ provides
\begin{equation}
 d_{j} = (\vg_{j}^{T}\mSigma^{-1}\vg_{j})^{-1} \quad j=1,\ldots,n
\end{equation}
As for $\mD_{2}$ and $\mGtilde_{22}$ the minimum is achieved when
\begin{equation}
\mGtilde_{22} \mD_{2} \mGtilde_{22}^{T} = \mL_{22}\mL_{22}^{T} = \mSigma_{2.1}
\end{equation}
Of course  these values are hypothetical and cannot be computed since they depend on $\mSigma$ but the FSOPT estimator will serve as a reference as the loss function is minimized for each realization of $\mX$.

\subsection{Oracle estimator}
The Oracle estimator is obtained by minimizing $\risk{\mGtilde \mD \mGtilde^{T},\mSigma}=\E{\loss{\mGtilde \mD \mGtilde^{T},\mSigma}}$. Since the columns of $\mX$ are independent and follow a multivariate Gaussian distribution with zero mean and covariance matrix $\mSigma$, which we denote as $\mX \dist \mCN{p,n}{\mzero}{\mSigma}{\eye{n}}$, it follows that $\mX \dist \mL \mXbar$ where $\mXbar \dist \mCN{p,n}{\mzero}{\eye{p}}{\eye{n}}$ where $\dist$ means ``has the same distribution as ''. Consequently $\mS \dist \mL \mWbar \mL^{T}$ where $\mWbar$ follows a singular Wishart distribution \cite{Srivastava03}. We let $\mWbar = \mGbar \mGbar^{T}$ with $\mGbar$ a lower triangular matrix with positive diagonal elements. Then all elements of $\mGbar$ are independent and $\mGbar_{ij} \dist \Gaussian{0}{1}$ for $i>j$, $\mGbar_{jj} \dist \sqrt{\chisquare{n-j+1}}$. Moreover we have $\mG \dist \mL\mGbar$.  It follows that the loss function in \eqref{loss_function} can be rewritten as
\begin{align}
\loss{\mGtilde \mD \mGtilde^{T},\mSigma} &= \tr{\mGbar \mD_{1} \mGbar^{T}} - \log \det{\mD_{1} } - \log \det{\mGbar_{11} \mGbar_{11}^{T}} \nonumber \\
&+ \tr{\mL_{22}^{-1}\mGtilde_{22} \mD_{2} \mGtilde_{22}^{T}\mL_{22}^{-T}} - \log \det{\mL_{22}^{-1}\mGtilde_{22} \mD_{2} \mGtilde_{22}^{T}\mL_{22}^{-T}}
\end{align}
Now
\begin{align}
\E{\tr{\mGbar \mD_{1} \mGbar^{T}}} &= \E{\sum_{i \geq j} d_{j} \mGbar_{ij}^{2}} \nonumber \\
&= \sum_{i > j} d_{j} + \sum_{j=1}^{n}d_{j}\E{\chisquare{n-j+1}} \nonumber \\
&= \sum_{j=1}^{n} (p-j+\E{\chisquare{n-j+1}} )d_{j} \\
\E{\log \det{\mGbar_{11} \mGbar_{11}^{T}}} &= \sum_{j=1}^{n} \E{\log \chisquare{n-j+1}}
\end{align}
which implies that
\begin{align}
\risk{\mGtilde \mD \mGtilde^{T},\mSigma} &=-\sum_{j=1}^{n} \E{\log \chisquare{n-j+1}} +\sum_{j=1}^{n} (p+n-2j+1)d_{j} - \sum_{j=1}^{n} \log d_{j} \nonumber \\
&+ \tr{\mL_{22}^{-1}\mGtilde_{22} \mD_{2} \mGtilde_{22}^{T}\mL_{22}^{-T}} - \log \det{\mL_{22}^{-1}\mGtilde_{22} \mD_{2} \mGtilde_{22}^{T}\mL_{22}^{-T}}
\end{align}
The risk is thus minimized for $d_{j}=(p+n-2j+1)^{-1}$, $j=1,\ldots,n$ and $\mGtilde_{22} \mD_{2} \mGtilde_{22}^{T} = \mL_{22}\mL_{22}^{T} $. In contrast to FSOPT, the Oracle estimator provides a  value of $\mD_{1}$ that can be computed. Yet, similarly to the FSOPT estimator, $\mGtilde_{22}$ and $\mD_{2}$ still depend on the unknown $\mSigma_{2.1}$. However this was expected since nothing can be inferred about $\mSigma_{2.1}$ with only $n<p$ observations.

\subsection{Practical scheme}
In order to come up with a bona fide estimator we first set $d_{j}=(p+n-2j+1)^{-1}$, $j=1,\ldots,n$ as in the Oracle estimator. Then we must decide upon $\mD_{2}$ and $\mGtilde_{22}$. Their optimal choice is $\mGtilde_{22} \mD_{2} \mGtilde_{22}^{T} = \mL_{22}\mL_{22}^{T} $ which cannot be met in practice. Therefore we must set these matrices arbitrarily, just as a target matrix is chosen in conventional linear shrinkage. Before proceeding, our intuition is that the least influential is $\mSigma_{2.1}$ the least the impact of a wrong guess for $\mD_{2}$ and $\mGtilde_{22}$. Our idea is to use a Cholesky decomposition of $\mSigma$ with complete pivoting \cite{Higham92} so that the diagonal elements of the so-obtained lower triangular factor are in decreasing order. Alternatively a QR decomposition with pivoting \cite{Golub96} of $\mX^{T}$ can be computed as $\mX^{T} \mPi = \mQ \mH^{T}$ where $\mPi$ is a $p \times p$ permutation matrix, $\mQ$ is a $n \times n$ orthogonal matrix and $\mH$ is a $p \times n$ lower triangular with positive and decreasing diagonal elements. This amounts to temporarily replacing $\mX$ by $\mY = \mPi^{T} \mX$, which corresponds to a permutation of the rows of $\mX$. The sample covariance matrix of $\mY$ is given by $\mY\mY^{T} = \mPi^{T} \mX \mX^{T} \mPi = \mH \mQ^{T} \mQ \mH^{T} = \mH \mH^{T}$ and hence $\mH$ is the Cholesky factor of $\mY\mY^{T}$ with decreasing diagonal elements. If we partition $\mH = \begin{pmatrix} \mH_{11}  \\ \mH_{21}  \end{pmatrix}$, then $\mH$ can be augmented to
\begin{equation}
\mHtilde = \begin{pmatrix} \mH_{11}  & \mzero \\ \mH_{21}  & \alpha\eye{p-n}  \end{pmatrix}
\end{equation}
where $\alpha = \mH_{11} (n,n)$ with the underlying  idea that the diagonal elements of $\mHtilde$ should decrease smoothly, so that we set the $p-n$ last diagonal elements to the minimal value of the $n$ first diagonal elements. In the same spirit we set $\mD_{2}=\beta \eye{p-n}$ with $\beta = d_{n}=(p-n+1)^{-1}$. The regularized estimate of the covariance matrix of $\mY$ is thus $\mHtilde \mD \mHtilde^{T}$ and thus, coming back to $\mX$, we obtain
\begin{equation}\label{Sigmahat}
\mSigmahat = \mPi \begin{pmatrix} \mH_{11}  & \mzero \\ \mH_{21}  & \alpha\eye{p-n}  \end{pmatrix}   \begin{pmatrix} \mD_{1}  & \mzero \\ \mzero  & \beta \eye{p-n}  \end{pmatrix} \begin{pmatrix} \mH_{11}  & \mH_{21}^{T} \\ \mzero  & \alpha\eye{p-n}  \end{pmatrix} \mPi^{T}
\end{equation}
Note that the regularizing parameters $\alpha$ and $\beta$ are chosen automatically: $\alpha$ will depend on the data while $\beta$ is fixed. The scheme allows to obtain a positive definite covariance matrix estimate at a reasonable cost since only a Cholesky decomposition is required. 

\section{Numerical simulations}
In this section we study the performance of $\mSigmahat$ in \eqref{Sigmahat} by evaluating its Stein's risk. Through preliminary simulations we observed that the FSOPT and the Oracle estimators perform the same so we only consider the latter. For comparison purposes we also evaluate Stein's risks of the benchmark linear shrinkage estimator of \cite{Ledoit04} (LW-LS in the figures), the nonlinear shrinkage estimator of \cite{Ledoit17} (LW-NLS) and the methods of \cite{Bickel08} (BL) and \cite{Rothman10} (RLZ). We consider a scenario where $p=200$. As for the spectrum of $\mSigma$ it is divided in two parts:  a set of  $\eta p$ large eigenvalues uniformly distributed on $[\lambdamax/2,\lambdamax]$ and a set of $(1-\eta) p$ small eigenvalues  uniformly distributed on $[0.5,1]$. $\eta$ controls the shape of the spectrum of $\mSigma$ while $\lambdamax$ controls its condition number, say $\cond{\mSigma}$. The influence of these parameters is seldom studied, yet we show that they have an important impact on which method should be retained.

\begin{figure}[p]
\centering
\subfigure{\includegraphics[width=12cm]{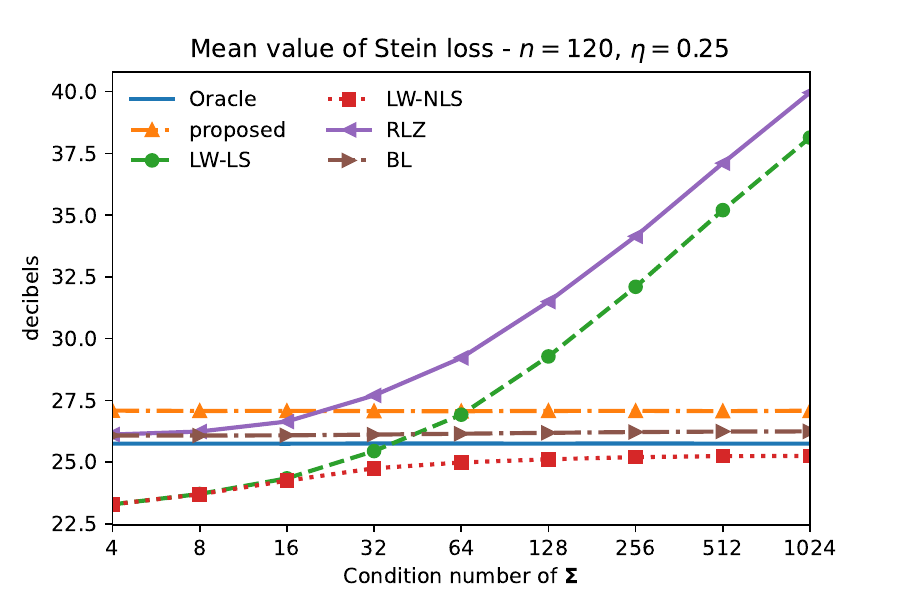}} \\
\subfigure{\includegraphics[width=12cm]{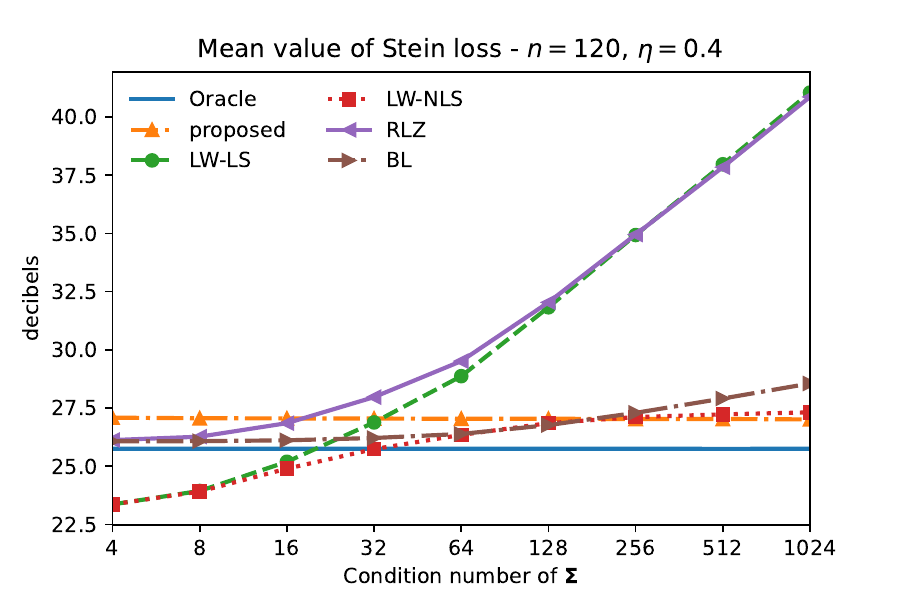}} \\
\caption{Mean value of Stein's loss versus $\cond{\mSigma}$. $n=120$ and varying $\eta$.}
\label{fig:risk_vs_condnum_n=120}
\end{figure}

In a first experiment we fix $n=120$, we vary the condition number between $4$ and $1024$ and we consider two different values of $\eta$, namely  $\eta=0.25$ and $\eta=0.4$. The results are reported in Figure \ref{fig:risk_vs_condnum_n=120}. Several key observations can be made. First, we observe that LS and RLZ are highly sensitive to the condition number of $\mSigma$ and also depends on the shape of the spectrum.  LW-LS performs better than $\mSigmahat$ only when $\cond{\mSigma}$ is below some threshold, and this threshold decreases as $\eta$ increases. In most situations the risk of $\mSigmahat$ is smaller than that of LW-LS.   In contrast, the risks of our method, BL and LW-NLS are independent of  $\cond{\mSigma}$ and of $\eta$, which is a very desirable feature in practice. Note that LW-NLS is the best estimator and our method performs slightly worse than BL, except with large condition numbers. However BL requires the user to set the regularized Cholesky factor bandwidth while our method is fully automatic and parameter-free.

\begin{figure}[p]
\centering
\subfigure{\includegraphics[width=12cm]{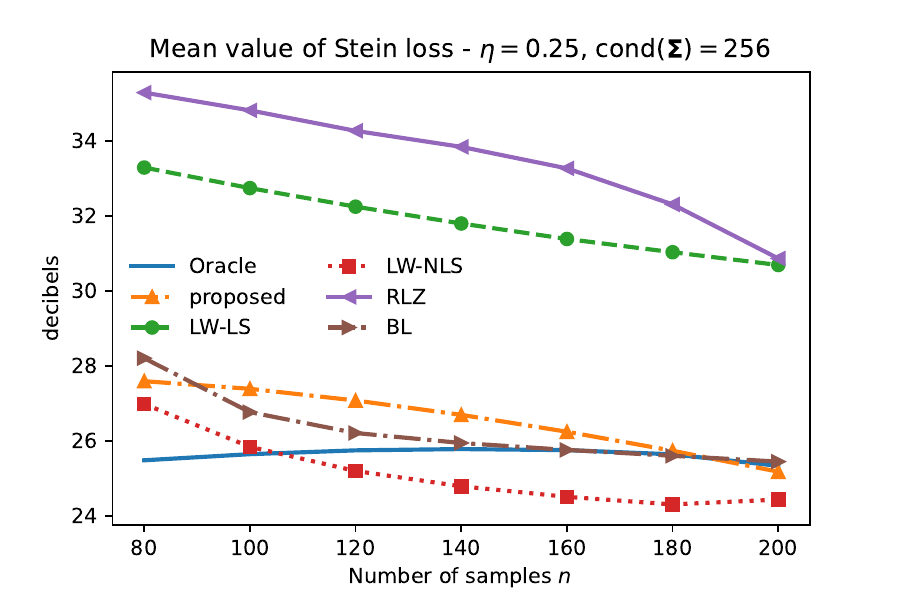}}	 \\
\subfigure{\includegraphics[width=12cm]{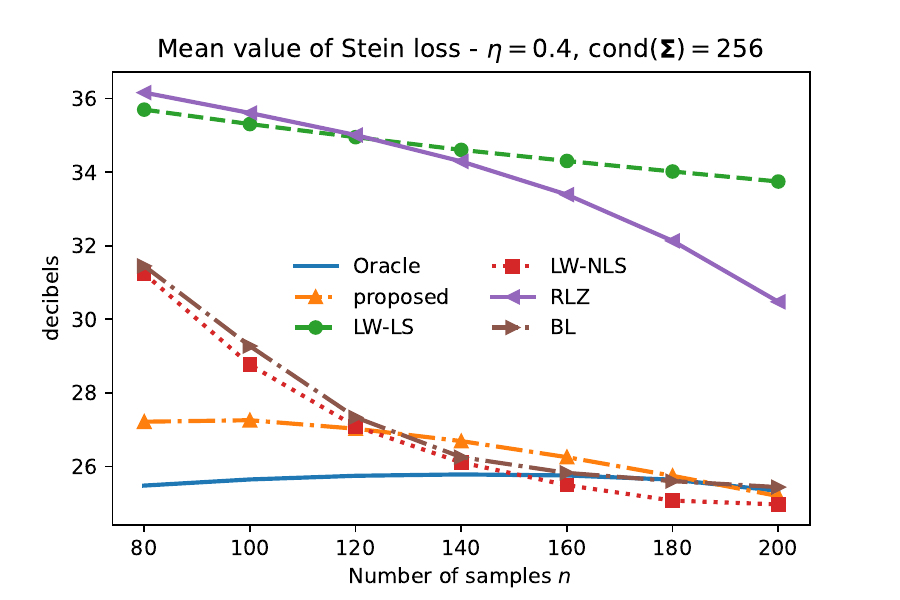}}
\caption{Mean value of Stein's loss versus $n$.  $\cond{\mSigma}=256$ and varying $\eta$.}	
\label{fig:risk_vs_n_condnum=256}
\end{figure}

In Figure \ref{fig:risk_vs_n_condnum=256} we study the influence of $n$ when the condition number of $\mSigma$ is held constant at $\cond{\mSigma}=256$. It can be seen that LW-LS and RLZ results in the largest risks while LW-NLS provides the smallest risk, except when $\eta=0.4$ and $n \leq 120$ where our method performs better. Again we observe a similar risk between our method and BL, unless $n$ is small and $\eta$ large.

\begin{figure}[p]
\centering
\includegraphics[width=12cm]{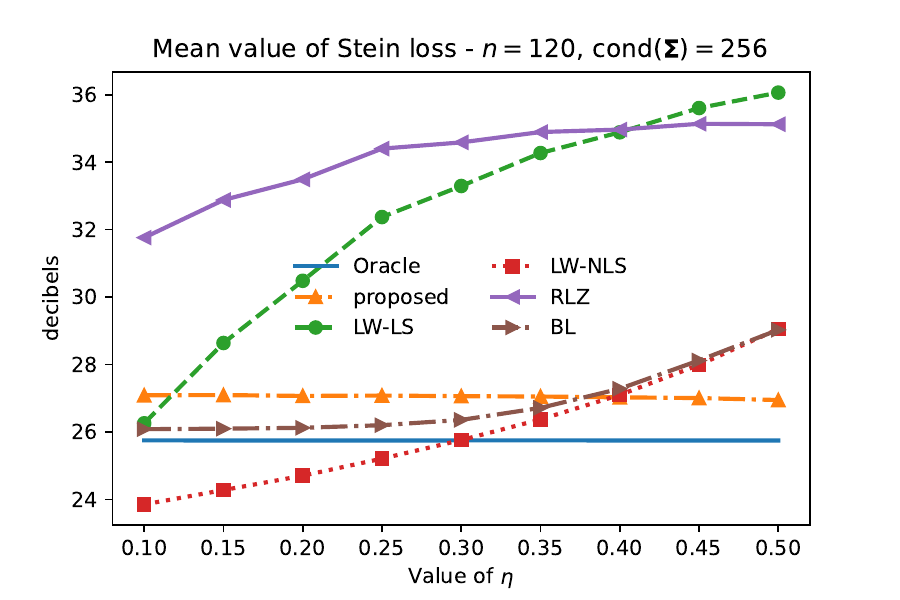}
\caption{Mean value of Stein's loss versus $\eta$.  $n=120$ and $\cond{\mSigma}=256$.}	
\label{fig:risk_vs_eta_n=120_condnum=256}
\end{figure}

Finally the influence of $\eta$ is investigated in Figure \ref{fig:risk_vs_eta_n=120_condnum=256}. It shows that LW-NLS performs very well for $\eta \leq 0.3$. Our method has a constant risk over $\eta$ and so it provides an interesting alternative for $\eta \geq 0.4$.

\section{Conclusions}
We addressed covariance matrix estimation in the singular case ($n < p$) through regularization of the rank-deficient Cholesky factor of the sample covariance matrix. We introduced the finite sample optimum estimator, which minimizes Stein's loss for every data matrix, the Oracle estimator, which minimizes the average loss, and a new practical scheme obtained by augmenting the Cholesky factor. The latter method results in a rather simple,  positive definite covariance matrix estimate whose regularization parameters are chosen automatically. We illustrated the impact of the condition number of the population covariance matrix as well as of the shape of its spectrum. We showed that the estimator proposed here provides an interesting solution when $n$ is small and   the large eigenvalues occupy at least $30\%$ of the spectrum.

\begin{thebibliography}{10}
	
	\bibitem{Pourahmadi13}
	M.~Pourahmadi.
	\newblock {\em High dimensional covariance estimation}.
	\newblock Wiley series in probability and statistics. John Wiley \& Sons,
	Hoboken, NJ, 2013.
	
	\bibitem{Zagidullina21}
	A.~Zagidullina.
	\newblock {\em High-dimensional covariance matrix estimation - {A}n
		introduction to random matrix theory}.
	\newblock SpringerBriefs in Applied Statistics and Econometrics, 2021.
	
	\bibitem{Ledoit22b}
	O.~Ledoit and M.~Wolf.
	\newblock The power of (non-)linear shrinking: {A} review and guide to
	covariance matrix estimation.
	\newblock {\em Journal of Financial Econometrics}, 20(1):187--218, 2022.
	
	\bibitem{Stein86}
	C.~Stein.
	\newblock Lectures on the theory of estimation of many parameters.
	\newblock {\em Journal of Mathematical Sciences}, 34:1373--1403, July 1986.
	
	\bibitem{Dey85}
	D.~K. Dey and C.~Srinivasan.
	\newblock {E}stimation of a covariance matrix under {S}tein's loss.
	\newblock {\em The Annals of Statistics}, 13(4):1581--1591, December 1985.
	
	\bibitem{Dey86}
	D.~K. Dey and C.~Srinivasan.
	\newblock Trimmed minimax estimator of a covariance matrix.
	\newblock {\em Annals Institute Statistical Mathematics}, 38:101--108, 1986.
	
	\bibitem{Sheena92}
	Y.~Sheena and A.~Takemura.
	\newblock Inadmissibility of non-order preserving orthogonally invariant
	estimators of the covariance matrix in the case of {S}tein's loss.
	\newblock {\em Journal of Multivariate Analysis}, 41:117--131, 1992.
	
	\bibitem{Perron92}
	F~Perron.
	\newblock Minimax estimators of a covariance matrix.
	\newblock {\em Journal of Multivariate Analysis}, 43(1):16 -- 28, 1992.
	
	\bibitem{Ledoit04}
	O.~Ledoit and M.~Wolf.
	\newblock A well-conditioned estimator for large-dimensional covariance
	matrices.
	\newblock {\em Journal of Multivariate Analysis}, 88(2):365--411, February
	2004.
	
	\bibitem{Ledoit12}
	O.~Ledoit and M.~Wolf.
	\newblock Nonlinear shrinkage estimation of large-dimensional covariance
	matrices.
	\newblock {\em The Annals of Statistics}, 40(2):1024--1060, April 2012.
	
	\bibitem{Ledoit17}
	O.~Ledoit and M.~Wolf.
	\newblock Direct nonlinear shrinkage estimation of large-dimensional covariance
	matrices.
	\newblock Technical report, University of Zurich, 2017.
	\newblock Working paper no. 264.
	
	\bibitem{Ledoit18}
	O.~Ledoit and M.~Wolf.
	\newblock Optimal estimation of a large dimensional covariance matrix under
	{S}tein's loss.
	\newblock {\em Bernoulli}, 24(4B):3791--3832, November 2018.
	
	\bibitem{Ledoit20}
	O.~Ledoit and M.~Wolf.
	\newblock Analytical nonlinear shrinkage of large-dimensional covariance
	matrices.
	\newblock {\em The Annals of Statistics}, 48(5):3043--3065, October 2020.
	
	\bibitem{Ledoit22}
	O.~Ledoit and M.~Wolf.
	\newblock Quadratic shrinkage for large covariance matrices.
	\newblock {\em Bernoulli}, 28(3):1519--1547, August 2022.
	
	\bibitem{Stein56}
	C.~Stein.
	\newblock Inadmissibility of the usual estimator for the mean of a multivariate
	distribution.
	\newblock In {\em Proceedings 3rd Berkeley Symposium on Mathematical Statistics
		and Probability}, pages 197--206, 1956.
	
	\bibitem{James61}
	W.~James and C.~Stein.
	\newblock Estimation with quadratic loss.
	\newblock In {\em Proceedings 4th Berkeley Symposium on Mathematical Statistics
		and Probability}, pages 361--380, 1961.
	
	\bibitem{Stein81}
	C.~Stein.
	\newblock Estimation of the mean of a multivariate normal distribution.
	\newblock {\em The Annals of Statistics}, 9(6):1135--1151, November 1981.
	
	\bibitem{Ledoit21}
	O.~Ledoit and M.~Wolf.
	\newblock Shrinkage estimation of large covariance matrices: Keep it simple,
	statistician?
	\newblock {\em Journal of Multivariate Analysis}, 186:104796, November 2021.
	
	\bibitem{Bodnar14}
	T.~Bodnar, A.~K. Gupta, and N.~Parolya.
	\newblock On the strong convergence of the optimal linear shrinkage estimator
	for large dimensional covariance matrix.
	\newblock {\em Journal of Multivariate Analysis}, 132:215--228, November 2014.
	
	\bibitem{Selliah64}
	J.~B. Selliah.
	\newblock Estimation and testing problems in a {W}ishart distribution.
	\newblock Technical report no. 10, Department of Statistics, Stanford
	University, January 10 1964.
	
	\bibitem{Eaton87}
	M.~L. Eaton and I.~Olkin.
	\newblock Best equivariant estimators of a {C}holesky decomposition.
	\newblock {\em The Annals of Statistics}, 15(4):1639--1650, 1987.
	
	\bibitem{Tsukuma16c}
	H.~Tsukuma and T.~Kubokawa.
	\newblock Unified improvements in estimation of a normal covariance matrix in
	high and low dimensions.
	\newblock {\em Journal of Multivariate Analysis}, 143:233--248, January 2016.
	
	\bibitem{Pourahmadi99}
	M.~Pourahmadi.
	\newblock Joint mean-covariance models with applications to longitudinal data:
	Unconstrained parameterisation.
	\newblock {\em Biometrika}, 86(3):677--690, September 1999.
	
	\bibitem{Pourahmadi00}
	M.~Pourahmadi.
	\newblock Maximum likelihood estimation of generalised linear models for
	multivariate normal covariance matrix.
	\newblock {\em Biometrika}, 87(2):425--435, June 2000.
	
	\bibitem{Wu03}
	W.~B. Wu and M.~Pourahmadi.
	\newblock Nonparametric estimation of large covariance matrices of longitudinal
	data.
	\newblock {\em Biometrika}, 90(4):831--844, December 2003.
	
	\bibitem{Huang06}
	J.~Z. Huang, N.~Liu, M.~Pourahmadi, and L.~Liu.
	\newblock Covariance matrix selection and estimation via penalised normal
	likelihood.
	\newblock {\em Biometrika}, 93(1):85--98, March 2006.
	
	\bibitem{Bickel08}
	P.~J. Bickel and E.~Levina.
	\newblock Regularization of large covariance matrices.
	\newblock {\em The Annals of Statsitics}, 36(1):199--227, February 2008.
	
	\bibitem{Rothman10}
	A.~J. Rothman, E.~Levina, and J.~Zhu.
	\newblock A new approach to {C}holesky-based covariance regularization in high
	dimensions.
	\newblock {\em Biometrika}, 97(3):539--550, September 2010.
	
	\bibitem{Srivastava03}
	M.~S. Srivastava.
	\newblock Singular {W}ishart and multivariate beta distributions.
	\newblock {\em The Annals of Statistics}, 31(5):1537--1560, October 2003.
	
	\bibitem{Higham92}
	N.~J. Higham.
	\newblock {\em Accuracy and stability of numerical algorithms}.
	\newblock Society for Industrial and Applied Mathematics, Philadelphia, 1992.
	
	\bibitem{Golub96}
	G.~Golub and C.~Van Loan.
	\newblock {\em {M}atrix {C}omputations}.
	\newblock {J}ohn {H}opkins {U}niversity {P}ress, {B}altimore, 3rd edition,
	1996.
	
\end{thebibliography}

\end{document}